\documentclass[a4paper,10pt]{article}
\usepackage[english]{babel}
\usepackage[utf8]{inputenc}

\usepackage{amsmath}
\usepackage{amsfonts}
\usepackage{amssymb}
\usepackage{amsthm}
\usepackage{graphicx}
\usepackage{pdfpages}

\newcommand{\mytext}[1]{ \: \textrm{#1} \: }
\newcommand{\mysetdescr}[2]{\left\{ { #1 \: \left| \: #2 \right. } \right\} }

\newcommand{\myfr}[1]{ \mathfrak{#1} }
\newcommand{\mydownarrow}{{\downarrow \,}}
\newcommand{\myuparrow}{{\uparrow \,}}
\newcommand{\myouparrow}{\uparrow_{_{_{\!\!\!\circ}}}}
\newcommand{\myodownarrow}{\downarrow^{^{\!\!\!\circ}}}

\newcommand{\myN}{\mathbb{N}}
\newcommand{\myNk}[1]{\underline {#1}}
\newcommand\mytimes{{\times}}

\newcommand\myurbild[1]{#1^{-1}}

\def\A{{\cal A}}
\def\C{{\cal C}}
\def\E{{\cal E}}
\def\H{{\cal H}}
\def\I{{\cal I}}
\def\S{{\cal S}}
\def\U{{\cal U}}

\newcommand{\mfP}{ \myfr{P} }

\newcommand{\mygxix}[2]{G_{#1}(#2)}
\newcommand{\gxix}{\mygxix{\xi}{x}}
\newcommand{\grxix}{\mygxix{\rho(\xi)}{x}}

\newcommand{\gzex}{\mygxix{\zeta}{x}}

\newcommand{\mydngxix}[2]{{\myodownarrow} \mygxix{#1}{#2}}
\newcommand{\myupgxix}[2]{{\myouparrow} \mygxix{#1}{#2}}

\newcommand{\dngxix}{{\myodownarrow} \gxix}
\newcommand{\upgxix}{{\myouparrow} \gxix}

\newcommand{\myapxix}[3]{\alpha_{#1, #2}(#3)}
\newcommand{\myapexix}[3]{\alpha_{#1, \eta(#2)}(#3)}

\newcommand{\myepxix}[3]{\eta_{#1}(#2)(#3)}

\newcommand{\myaxix}[2]{\alpha_{#1}\left(#2\right)}
\newcommand{\myarxix}[2]{\alpha_{\rho(#1)}(#2)}

\newcommand{\axix}{\myaxix{\xi}{x}}

\newcommand{\axiy}{\myaxix{\xi}{y}}

\newcommand{\eps}{\epsilon}

\newcommand{\mf}[1]{\mathfrak{ #1 }}
\newcommand{\fa}{\mf{a}}
\newcommand{\fb}{\mf{b}}

\newcommand{\exix}{\eta(\xi)(x)}

\newcommand{\arie}{P \in \mfP_r$, $\xi \in \H(P,R)$, $x \in P}

\newcommand{\spc}[1]{{\; #1 \;}}

\def\BP{\begin{proof}}
\def\EP{\end{proof}}

\DeclareMathOperator{\lcm}{lcm}

\begin{document}

\theoremstyle{plain}
\theoremstyle{plain}
\newtheorem{condition}{Condition}
\newtheorem{theorem}{Theorem}
\newtheorem{definition}{Definition}
\newtheorem{corollary}{Corollary}
\newtheorem{lemma}{Lemma}
\newtheorem{proposition}{Proposition}

\title{\bf Calculation Rules and Cancellation Rules for Strong Hom-Schemes}

\author{\sc Frank a Campo}
\date{\small Seilerwall 33, D 41747 Viersen, Germany\\
{\sf acampo.frank@gmail.com}}

\maketitle

\begin{abstract}
\noindent Let $\H(A,B)$ denote the set of homomorphisms from the poset $A$ to the poset $B$. In previous studies, the author has started to analyze what it is in the structure of finite posets $R$ and $S$ that results in $\# \H(P,R) \leq \# \H(P,S)$ for every finite poset $P$, if additional regularity conditions are imposed. In the present paper, it is examined if this relation (with or without regularity conditions) is compatible with the operations of order arithmetic and if cancellation rules hold.
\newline

{\bf Key words:} 
\end{abstract}

\section{Introduction} \label{sec_introduction}

\begin{figure} 
\begin{center}
\includegraphics[trim = 70 710 200 70, clip]{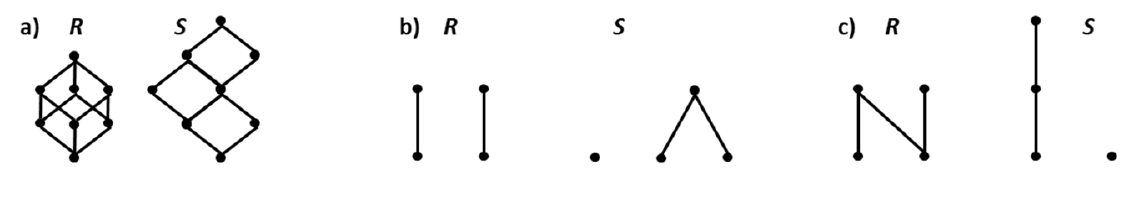}
\caption{\label{figure_1} Three examples for posets $R$ and $S$ with $\# \H(P,R) \leq \# \H(P,S)$ for every finite poset $P$ \cite{aCampo_toappear_1}.}
\end{center}
\end{figure}

Based on a theorem of Lov\'{a}sz \cite{Lovasz_1967} and an own observation \cite[Theorem 5]{aCampo_2018}, the author \cite{aCampo_toappear_1,aCampo_toappear_2} has worked about homomorphism sets under the aspect:{\em  What is it in the structure of finite poests $R$ and $S$ that results in $\# \H(P,R) \leq \# \H(P,S)$ for every finite poset $P$?} ($\H(A,B)$ denotes the set of homomorphisms from the poset $A$ to the poset $B$.) Three nontrivial examples for such pairs of posets are shown in Figure \ref{figure_1}. More can be found in \cite{aCampo_toappear_1}, where a detailed survey over the research about homomorphism sets and its results is also presented.

The systematic examination of the phenomenon ``$\# \H(P,R) \leq \# \H(P,S)$ for every finite poset $P$'' in \cite{aCampo_toappear_1,aCampo_toappear_2} is based on the concept of the {\em strong Hom-scheme from $R$ to $S$}. In the case of existence, such a Hom-scheme defines a one-to-one mapping $\rho_P : \H(P,R) \rightarrow \H(P,S)$ for every finite poset $P \in \mfP_r$, where $\mfP_r$ is a representation system of the non-isomorphic finite posets. Postulating regularity conditions for the way, how $\rho_P$ maps homomorphimss from $\H(P,R)$ to $\H(P,S)$, gives rise to the definition of the {\em strong G-scheme} and the {\em strong I-scheme}. The relations $R \sqsubseteq S$ (there exists a strong Hom-scheme from $R$ to $S$), $R \sqsubseteq_G S$ (existence of a strong G-scheme), and $R \sqsubseteq_I S$ (existence of a strong I-scheme), turn out to be partial order relations on $\mfP_r$.

The present paper deals with calculation and cancellation rules for strong Hom-, G-, and I-schemes. Let ${\preceq} \in \{ \sqsubseteq, \sqsubseteq_G, \sqsubseteq_I  \}$. Our questions are:
\begin{itemize}
\item For finite posets $R, S$ with $R \preceq S$, what is the relation between the duals of $R$ and $S$, and how is $R \preceq S$ connected with $\H(Q,R) \preceq \H(Q,S)$ for every $Q \in \mfP$?
\item Let $R_1, R_2, S_1, S_2$ be finite posets with $R_1 \preceq S_1$ and $R_2 \preceq S_2$. Does this imply $R_1 \odot R_2 \preceq S_1 \odot S_2$, where $ \odot $ denotes the poset-operator direct sum, ordinal sum, or product?
\item Let $Q, R, S$ be finite posets with $ Q \odot R \preceq Q \odot S$, $\odot$ as above. When does $R \preceq S$ hold?
\end{itemize}

After the preparatory Sections \ref{sec_notation} and \ref{sec_connectivity}, the required apparatus about Hom-schemes is recalled from \cite{aCampo_toappear_1,aCampo_toappear_2} in Section \ref{sec_EV_HomS}. Additionally, it is shown in Proposition \ref{prop_zshP_reicht} that a (strong) Hom-scheme, G-scheme, or I-scheme only defined for connected finite posets can be extended to all finite posets. 

Section \ref{sec_calcrules} containes the calculation rules. They work without exceptions for strong Hom-schemes and strong G-schemes, but for strong I-schemes, we did not prove a calculation rule for ordinal sums and products.

The Sections \ref{sec_cancel_dirsum}, \ref{sec_cancel_ordsum}, and \ref{sec_cancel_product} contain the cancellation rules. For the direct sum of posets (Section \ref{sec_cancel_dirsum}), the cancellation rules work with no ifs or buts, but for ordinal sums (Section \ref{sec_cancel_ordsum}) and products (Section \ref{sec_cancel_product}), we need additional assumptions for strong I-schemes and partly also for strong G-schemes. Most of the proofs are quite technical, but the ansatz is similar in many of them, as explained in Section \ref{sec_ansatz_cancel}. However, there remains a gap: We did not succeed in proving a cancellation rule for $Q \oplus R \sqsubseteq Q \oplus S$.

\section{Basics and Notation} \label{sec_notation}

For a given set $X$, a subset of $X \times X$ is called a {\em (binary) relation} $R$ on $X$ and $X$ is called the {\em carrier} of $R$. The {\em induced relation} on $A \subseteq X$ is $R \cap ( A \mytimes A )$. We write $x R y$ for $(x,y) \in R$. The {\em dual relation} $R^d$ is defined by $x R^d y \Leftrightarrow y R x$ for all $x, y \in X$. A reflexive, antisymmetrical, and transitive relation $R$ on $X$ is called a {\em partial order relation}, and the pair $P = (X,R)$ is called a {\em partially ordered set} or simply a {\em poset}. For an reflexive and antisymmetrical relation (hence, in particular, for a partial order relation), we use the symbol $\leq$. ``$x < y$'' means ``$x \leq y$ and $x \not= y$''.

A subset $A \subseteq X$ is called a {\em chain} iff $x R y$ or $y R x$ for all $x, y \in A$; the cardinality of a chain is called its {\em length}. If a relation $R$ contains a finite longest chain, the length $h_R$ of such a chain is called the {\em height} of $R$.

For relations $R$ and $S$ with disjoint carriers $X$ and $Y$, their {\em direct sum} $R + S$ is defined as $R \cup S$, and their {\em ordinal sum} as $R \oplus S \equiv R \cup S \cup (X \times Y)$. Given posets $P = (X, \leq_P ) $ and $Q = (Y, \leq_Q)$, their product $P \mytimes Q \equiv (X \mytimes Y, \leq_{P \mytimes Q}) $ on $X \mytimes Y$ is defined by $(x_1,y_1) \leq_{P \mytimes Q} (x_2, y_2) $ iff $ x_1 \leq_P x_2$ and $y_1 \leq_Q y_2$.

For relations $R$ and $S$ on $X$ and $Y$, respectively, a mapping $\xi : X \rightarrow Y$ is called a {\em homomorphism}, iff $\xi(x) S \xi(y)$ holds for all $x, y \in X$ with $x R y$. A homomorphism $\xi$ is called {\em strict} iff it additionally fulfills $( x R y$ and $x \not= y ) \Rightarrow ( \xi(x) S \xi(y)$ and $\xi(x) \not= \xi(y) )$ for all $x, y \in X$, and it is called an {\em embedding} iff $\xi(x) S \xi(y) \Rightarrow x R y$ for all $x, y \in X$. Finally, an embedding is called an {\em isomorphism} iff it is onto. For posets $P$ and $Q$, we let $\H(P,Q)$ and $\S(P,Q)$ denote the set of homomorphims from $P$ to $Q$ and the set of strict homomorphims from $P$ to $Q$, respectively. $P \simeq Q$ indicates isomorphism. We equip homomorphism sets with the ordinary pointwise partial order.

$\mfP$ is the class of all {\em non-empty} finite posets, and $\mfP_r$ is a representation system of the non-isomorphic posets in $\mfP$. (Without loss of mathematical substance, we avoid repeated trivial distinctions of cases by excluding the empty poset.)

For a poset $P$ with carrier $X$, we use the notation $x \in P$ instead of $x \in X$, and for posets $P$ and $Q$ with carriers $X$ and $Y$, respectively, we write $\xi : P \rightarrow Q$ instead of $\xi : X \rightarrow Y$ for a homomorphism $\xi \in \H(P,Q)$.

Given a relation $R$ on $X$, we define for $A \subseteq X$
\begin{align*}
\mydownarrow A & \equiv \mysetdescr{ y \in X }{ \exists \; a \in A \; \mytext{:} \; y R a }, \\
{\myodownarrow} A & \equiv ( \mydownarrow A ) \setminus A, \\
\myuparrow A & \equiv \mysetdescr{ y \in X }{ \exists \; a \in A \; \mytext{:} \;  a R y }, \\
{\myouparrow} A & \equiv ( \myuparrow A ) \setminus A.
\end{align*}
For $x \in X$, we write $\mydownarrow x$ and ${\myodownarrow} x$ instead of $\mydownarrow \{x\}$ and ${\myodownarrow} \{x\}$, respectively, and correspondingly $\myuparrow x$ and ${\myouparrow} x$. If required, we label the arrows with the relation they are referring to. $A \subseteq X$ is called an {\em upset} iff $A = \myuparrow A$.

Additionally, we use the following notation from set theory:
\begin{align*}
\myNk{0} & \equiv  \emptyset \\
\myNk{n} & \equiv  \{ 1, \ldots, n \} \mytext{for every} n \in \myN.
\end{align*}

We denote the $i$\textsuperscript{th} component of $x \in X_1 \times \cdots \times X_I$ as usual by $x_i$. Due to $\H( P, Q \times R ) \simeq \H( P, Q ) \times \H(P, R)$ for all posets $P, Q, R$, we frequently write $( \xi_1, \xi_2)$ for a homomorphism $\xi \in \H( P, Q \times R )$. In order to avoid confusion in Section \ref{sec_cancel_product} (where we are dealing with triplets consisting of pairs and sets of pairs), we use additionally the canonical projections $\pi_1$ and $\pi_2$: for every $( a, b ) \in A \times B$, $\pi_1(a,b) \equiv a$ and $\pi_2(a,b) \equiv b$.

For sets $X$ and $Y$, $\A(X,Y)$ is the set of mappings from $X$ to $Y$. $id_X \in \A(X,X)$ is the {\em identity mapping}. For a mapping $f \in \A(X,Y)$, our symbols for the {\em pre-image} of $B \subseteq Y$ and of $y \in Y$ are
\begin{align*}
\myurbild{f}(B) & \equiv \mysetdescr{ x \in X }{ f(x) \in B}, \\
\myurbild{f}(y) & \equiv \myurbild{f} (\{ y\} ).
\end{align*}
For $A \subseteq X$ and $B \subseteq Y$ with $f(X) \subseteq B$, let $f \vert_A : A \rightarrow Y$ and $f \vert^B : X \rightarrow B$ denote the {\em pre-restriction} and {\em post-restriction} of $f$, respectively.

Let $A, B, C, D$ be sets with $A \cap C = \emptyset$. For mappings $f : A \rightarrow B$, $g : C \rightarrow D$, the mapping $ f \cup g : A \cup C \rightarrow B \cup D$ is defined by
\begin{align*}
( f \cup g )(x) & \equiv
\begin{cases}
f(x), & \mytext{if} x \in A; \\
g(x), & \mytext{if} x \in B.
\end{cases}
\end{align*}
Let $A, B$ be sets with $A \subseteq B$. For $f \in \A(A,B)$, we define recursively for every $a \in A$
\begin{align*}
f^0(a) & \equiv a, \\
\forall \; i \in \myN \; : \; f^{i+1}(a) & \equiv f( f^i( a ) ), \; \mytext{if} \; f^i( a ) \in A.
\end{align*}

Let $\I$ be a non-empty set, and let $N_i$ be a non-empty set for every $i \in \I$. The {\em Cartesian product} of the sets $N_i, i \in \I$, is defined as
\begin{eqnarray*}
\prod_{i \in \I} N_i & \; \equiv \; & 
\mysetdescr{ f \in \A \big( \I, \bigcup_{i \in \I} N_i \big)}{ f(i) \in N_i \mytext{for all} i \in \I }.
\end{eqnarray*}

\section{Connectivity} \label{sec_connectivity}

\begin{definition} \label{def_connected}
Let $P \in \mfP$, $A \subseteq P$, and $x, y \in A $. We say that $x$ and $y$ are {\em connected in $A$}, iff there are $z_0, z_1, \ldots , z_L \in A$, $L \in \myN_0$, with $x = z_0$, $y = z_L$ and $ z_{\ell-1} < z_\ell$ or $z_{\ell-1} > z_\ell$ for all $\ell \in \myNk{L}$. We call $z_0, \ldots , z_L$ a {\em zigzag line connecting $x$ and $y$}. We define for all $A \subseteq P$, $x \in A $, $B \subseteq A$
\begin{align*}
\gamma_A(x) & \equiv \mysetdescr{ y \in A }{ x \mytext{and} y \mytext{are connected in} A}, \\
\gamma_A(B) & \equiv \bigcup_{ b \in B } \gamma_A(b).
\end{align*}
\end{definition}
The sets $\gamma_P(x), x \in P$, are called the {\em connectivity components} of $P$. Every poset  is the direct sum of its connectivity components. A poset $P \in \mfP$ is {\em connected} iff $\gamma_P(x) = P$ for an/all $x \in P$, and a non-empty subset $A \subseteq P$ is called {\em connected (in $P$)} iff $\gamma_A(x) = A$ for an/all $x \in A$. For posets $P$ and $Q$, $P$ connected, the image $\xi(P)$ of $P$ under a homomorphism $\xi \in \H(P,Q)$ is connected in $Q$; in particular, $\xi(P)$ is a subset of a single connectivity component of $Q$. In consequence, if $P$ is a connected poset, then for all disjoint posets $R$ and $S$
\begin{align*}
\H(P, R + S) & \simeq \H(P,R) + \H(P,S),
\end{align*}
and correspondingly also for $\S(P, R + S )$. Furthermore, for posets $P$ and $Q$, every homomorphism $\xi \in \H(P,Q)$ can be decomposed into its restrictions to the connectivity components $K_i$ of $P$, $i \in \I$:
\begin{align*}
\xi & = \bigcup_{i \in \I} \xi \vert_{K_i}.
\end{align*}

The following results have been proven by the author \cite{aCampo_2018,aCampo_toappear_1}:
\begin{lemma} \label{lemma_ueber_gamma}
Let $A \subseteq P$. The relation ``connected in $A$'' is an equivalence relation on $A$ with partition $\mysetdescr{\gamma_A(a)}{a \in A}$. For $B \subseteq A \subseteq A' \subseteq P$ we have
\begin{align}
\gamma_A(B) & \subseteq \gamma_{A'}(B), \label{gamma_ABAB} \\
\gamma_A(B) & = \gamma_{\gamma_A(B)}(B). \label{gamma_ggc}
\end{align}
\end{lemma}

The next definition is fundamental also for the recent article:

\begin{definition}
Let $P = (X, \leq)$ be a finite poset, let $Y$ be a set, and let $\xi \in \A(X,Y)$ be a mapping. We define for all $x \in P$
\begin{align*}
\gxix & \equiv \gamma_{\myurbild{\xi}( \xi(x) )}( x ).
\end{align*}
\end{definition}
The set $\gxix$ is always a connected subset of $P$ with $x \in \gxix$. We have
\begin{align} \label{gxix_gammagxix}
\gxix & = \gamma_{\gxix}(x),
\end{align}
and $\xi$ is strict iff $\gxix = \{x\}$ for all $x \in P$. For a one-to-one mapping $\sigma : Y \rightarrow Z$, we have $\mygxix{\sigma \circ \xi}{x} = \gxix$ for every $\xi \in \A(X,Y), x \in P$.

We need two additional lemmata in this paper:

\begin{lemma} \label{lemma_G_prod}
Let $P \in \mfP$, and let $Y$ and $Z$ be non-empty finite sets. For $\xi \in \A(P, Y \mytimes Z)$ let $\xi_1 \in \A(P,Y)$ and $\xi_2 \in \A(P,Z)$ be defined by $( \xi_1(x), \xi_2(x) ) = \xi(x)$ for all $x \in X$. Then, for all $x \in X$,
\begin{align} \label{gxix_produkt}
\gxix & = \gamma_{ \mygxix{\xi_1}{x} \cap \mygxix{\xi_2}{x} }(x).
\end{align}
\end{lemma}

\BP For every $ (a,b) \in Y \mytimes Z $, we have $\myurbild{\xi}(a,b) = \myurbild{\xi_1}(a) \cap \myurbild{\xi_2}(b)$, which yields $\myurbild{\xi}( \xi(x) ) \subseteq \myurbild{\xi_j}(\xi_j(x))$ for every $x \in P$ and $j \in \myNk{2}$. With \eqref{gamma_ABAB}, we conclude $\gxix \subseteq \mygxix{\xi_1}{x} \cap \mygxix{\xi_2}{x}$, and we get
\begin{equation*}
\gxix \; \stackrel{\eqref{gamma_ggc}}{=} \; \gamma_{\gxix}(x) \;  \stackrel{\eqref{gamma_ABAB}}{\subseteq} \; \gamma_{ \mygxix{\xi_1}{x} \cap \mygxix{\xi_2}{x} }(x).
\end{equation*}
Now let $y \in \gamma_{ \mygxix{\xi_1}{x} \cap \mygxix{\xi_2}{x} }(x)$, and let $z_0, \ldots , z_L$ be a zigzag line connecting $x$ and $y$ in $\mygxix{\xi_1}{x} \cap \mygxix{\xi_2}{x}$. Then $\xi_1(z_\ell) = \xi_1(x)$ and $\xi_2(z_\ell) = \xi_2(x)$ for every $\ell \in \myNk{L} \cup \{ 0 \}$. The zigzag line $z_0, \ldots , z_L$ connects thus $x$ and $y$ in $\myurbild{\xi_1}(\xi_1(x)) \cap \myurbild{\xi_2}(\xi_2(x)) = \myurbild{\xi}(\xi(x))$, and we conclude $y \in \gxix$.

\EP
\begin{lemma} \label{gxix_neu}
For every $\arie$
\begin{align} \label{gxix_postrestr_oben}
\forall \; B \subseteq R \; : \; \xi(P) \subseteq B \quad & \Rightarrow \quad \mygxix{ \xi \vert^B }{x} = \gxix, \\
\label{gxix_postrestr_unten}
\forall \; A \subseteq P \; : \; \gxix \subseteq A \quad & \Rightarrow \quad \mygxix{ \xi \vert_A }{x} = \gxix,
\end{align}
in the partial order induced on $A$.
\end{lemma}
\BP Due to $\xi(P) \subseteq B$, we have $\myurbild{\left( {\xi \vert^B} \right)}( \xi \vert^B (x)) = \myurbild{\xi}( \xi (x))$, and \eqref{gxix_postrestr_oben} follows. In the case $\gxix \subseteq A$ we have $x \in A$, and $ \xi \vert_A (x)$ is well-defined. Now
\begin{align*}
\gxix & \subseteq A \cap \myurbild{\xi}( \xi(x) ) = 
\myurbild{\left( \xi \vert_A \right) }( \xi(x) ) = 
\myurbild{\left( \xi \vert_A \right) }( \xi \vert_A(x) ),
\end{align*}
thus $\gxix \stackrel{\eqref{gxix_gammagxix}}{=} \gamma_{\gxix}(x) \stackrel{\eqref{gamma_ABAB}}{\subseteq} \gamma_{\myurbild{\left( \xi \vert_A \right) }( \xi \vert_A(x) )}(x) = \mygxix{\xi \vert_A}{x}$. Furthermore,
\begin{align*}
\mygxix{\xi \vert_A}{x} & = 
\gamma_{ \myurbild{ \left( \xi \vert_A \right) }( \xi \vert_A ( x ) )}(x) =
\gamma_{ \myurbild{ \left( \xi \vert_A \right)  }( \xi( x ) )}(x) 
\stackrel{\eqref{gamma_ABAB}}{\subseteq}
\gamma_{ \myurbild{ \xi }( \xi( x ))}(x) = \gxix.
\end{align*}

\EP

\section{Hom-Schemes} \label{sec_EV_HomS}

\begin{figure} 
\begin{center}
\includegraphics[trim = 75 640 185 65, clip]{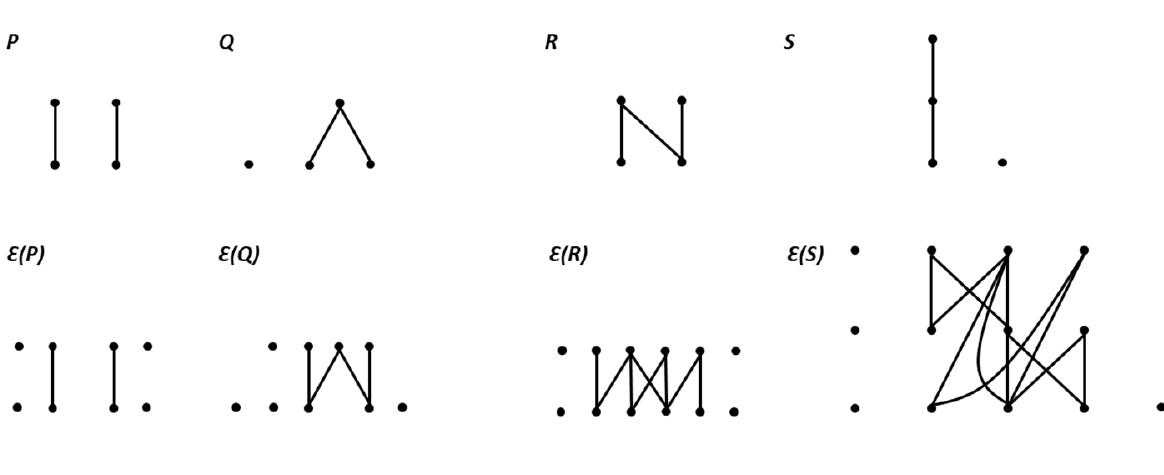}
\caption{\label{figure_2} Examples for EV-systems of posets \cite{aCampo_toappear_1}.}
\end{center}
\end{figure}

In in this section, the required definitions and results about Hom-schemes contained in \cite{aCampo_toappear_1,aCampo_toappear_2} are summarized without proofs. In Proposition \ref{prop_zshP_reicht}, a new result is presented.

\begin{definition} \label{def_EVsys}
Let $P$ be a poset. The {\em EV-system} $\E(P)$ of $P$ is defined as
\begin{align*}
\E(P) & \equiv \mysetdescr{ ( x, D, U ) }{ x \in P, D \subseteq  {\myodownarrow} x, U \subseteq {\myouparrow} x }.
\end{align*}
We equip $\E(P)$ with a reflexive and antisymmetrical relation: For all $\fa, \fb \in \E(P)$ we define
\begin{equation*}
\fa <_+ \fb \quad \equiv \quad \fa_1 \in \fb_2 \; \mytext{and} \; \fb_1 \in \fa_3,
\end{equation*}
and $\leq_+ \; \equiv \; <_+ \cup \; \mysetdescr{(\fa,\fa)}{\fa \in \E(P)}$.
\end{definition}

Four examples of EV-systems are shown in Figure \ref{figure_2}. ``EV'' reminds of the exploded view drawings in mechanical engineering: just as the EV-system $\E(P)$ does with the points of a poset $P$, the exploded view drawing of an engine shows the relationships between its components by distributing them in the drawing area in a well-arranged and meaningful way. 

The mapping $\E(P) \rightarrow P$ defined by $\fa \mapsto \fa_1$ is a strict homomorphism, and the mapping $P \rightarrow \E(P)$ with $x \mapsto \left( x, {\myodownarrow} x, {\myouparrow} x \right)$ is an embedding. We conclude $h_{\E(P)} = h_P$.

The following equations about the EV-systems of posets combined by order arithmetic follow directly from the definition. For every poset $P \in \mfP$ we have
\begin{align}
\E(P^d) & \simeq \E(P)^d, \label{EVsys_dual}
\end{align}
and for disjoint posets $P, Q \in \mfP$ we have
\begin{align}
\E(P+Q) = & \;\;\; \E(P) + \E(Q), \label{EVsys_dirsum} \\
\E( P \oplus Q ) = & \;\;\; \mysetdescr{ (x, D, U \cup V) }{(x,D,U) \in \E(P), V \subseteq Q} \label{EVsys_ordsum} \\
& \cup  \mysetdescr{ (x, D \cup E, U) }{(x,D,U) \in \E(Q), E \subseteq P}. \nonumber
\end{align}
Furthermore, $\E( P \times Q )$ is for any posets $P, Q \in \mfP$ the set of the $( (x,y), D, U )$ with
\begin{align}
(x,y) & \in P \times Q, \label{EVsys_prod} \\
D & \subseteq \left( ( \mydownarrow_P x ) \times ( \mydownarrow_Q y ) \right) \setminus \{ (x,y) \}, \nonumber \\
U & \subseteq \left( ( \myuparrow_P x ) \times ( \myuparrow_Q y ) \right) \setminus \{ (x,y) \}. \nonumber
\end{align}

Let $P, Q \in \mfP$ and $\xi \in \H(P,Q)$. We define a homomorphism $\alpha_{P,\xi} : P \rightarrow \E(Q)$ be setting for every $x \in P$
\begin{align*}
\myapxix{P}{\xi}{x} & \equiv \left( \xi(x), \xi( \dngxix ), 
                          \xi( \upgxix ) \right).
\end{align*}
If $P$ is fixed, we write $\axix$ instead of $\myapxix{P}{\xi}{x}$.

\begin{lemma}[{\cite[Proposition 2]{aCampo_toappear_1}}] \label{lemma_exix}
Let $R, S \in \mfP$, and let $\eps : \E(R) \rightarrow \E(S)$ be a strict homomorphism. For given $P \in \mfP$ and $\xi \in \H(P,R)$, we define for every $x \in P$
\begin{align*}
\eta(x) & \equiv \eps( \axix )_1,
\end{align*}
Then $\eta \in \H(P,S)$, and for every $x \in P$ we have
\begin{align*}
\mygxix{\eta}{x} & = \gxix.
\end{align*}
\end{lemma}

%\begin{definition} \label{def_Ma_ia}
%For $P \in \mfP$ with carrier $X$ we define for every $\fa \in \E(P)$
%\begin{align*}
%M(\fa) & \equiv \fa_2 \oplus \{ \fa_1 \} \oplus \fa_3,
%\end{align*}
%where $\fa_2$ and $\fa_3$ are treated as antichains. $\iota(\fa) : M(\fa) \rightarrow X $ is the canonical inclusion mapping of $M(\fa)$ in $X$.
%\end{definition}
%$M(\fa)$ is isomorphic to $A_{\# \fa_2} \oplus A_1 \oplus A_{\# \fa_3}$, and $\iota(\fa) \in \H(M(\fa),P)$ is a one-to-one homomorphism. We have $\fa = \myamaiax{\fa}{\fa_1}$ for every $\fa \in \E(P)$.

The central definition is:

\begin{definition} \label{def_RS_scheme}
Let $R, S \in \mfP$. We call a mapping
\begin{align*}
\rho & \in \prod_{P \in \mfP_r} \A( \H(P,R), \H(P,S) ) \quad \mytext{(Cartesian product)}
\end{align*}
a {\em Hom-scheme from $R$ to $S$}. We call a Hom-scheme $\rho$ from $R$ to $S$
\begin{itemize}
\item {\em strong} iff the mapping $\rho_P : \H(P,R) \rightarrow \H(P,S)$ is one-to-one for every $P \in \mfP_r$;
\item a {\em G-scheme} iff for every $P \in \mfP_r, \xi \in \H(P,R), x \in P$
\begin{align*}
G_{\rho_P(\xi)}(x) & = \gxix;
\end{align*}
\item {\em image-controlled} or an {\em I-scheme}, iff for every $P, Q \in \mfP_r$, $\xi \in \H(P,R), \zeta \in \H(Q,R)$, $x \in P$, $y \in Q$
\begin{equation*}
\myapxix{P}{\xi}{x} \leq_+ \myapxix{Q}{\zeta}{y} \; \Rightarrow \;
\myapxix{P}{\rho_P(\xi)}{x} \leq_+ \myapxix{Q}{\rho_Q(\zeta)}{y},
\end{equation*}
where ``$<_+$'' on the left side implies ``$<_+$'' on the right side.
\end{itemize}
If $P$ is fixed, we write $\rho(\xi)$ instead of $\rho_P(\xi)$.\end{definition}
If a one-to-one homomorphism $\sigma : R \rightarrow S$ exists, we get a strong I-scheme from $R$ to $S$ by defining $\rho_P(\xi) \equiv \sigma \circ \xi$ for all $P \in \mfP$, $ \xi \in \H(P,R)$. For an image-controlled Hom-scheme $\rho$, we have due to the antisymmetry of $\leq_+$
\begin{equation} \label{eq_imagebased_rhowert}
\myapxix{P}{\xi}{x} = \myapxix{Q}{\zeta}{y} 
\; \Rightarrow \;
\myapxix{P}{\rho_P(\xi)}{x} = \myapxix{Q}{\rho_Q(\zeta)}{y}. 
\end{equation}

An I-scheme is always a G-scheme \cite[Theorem 4]{aCampo_toappear_1}. We write 
\begin{equation*}
R \sqsubseteq S \; / \; R \sqsubseteq_G S \; / \; R \sqsubseteq_I S
\end{equation*}
iff a strong Hom-scheme / a strong G-scheme / a strong I-scheme exists from $R$ to $S$. The three relations $\sqsubseteq$, $\sqsubseteq_G$, and $\sqsubseteq_I$ define partial orders on $\mfP_r$ \cite[Theorems 2 and 3]{aCampo_toappear_1}.

The main results of \cite{aCampo_toappear_1,aCampo_toappear_2} about strong I-schemes and strong G-schemes are:

\begin{theorem}[{\cite[Theorems 5 and 7]{aCampo_toappear_1}}] \label{theo_strongI_equiv_eps}
For posets $R, S \in \mfP$, we have $ R \sqsubseteq_I S$ iff there exists a one-to-one homomorphism $\eps : \E(R) \rightarrow \E(S)$ with
\begin{equation} \label{eq_cond_univ_aexix}
\myapexix{P}{\xi}{x} \quad = \quad \eps( \myapxix{P}{\xi}{x}  ).
\end{equation}
for every $\arie$, where
\begin{equation*}
\exix \; \equiv \; \myepxix{P}{\xi}{x} \; \equiv \; \eps( \myapxix{P}{\xi}{x} )_1
\end{equation*}
for every $\arie$. In particular, $\eta$ is a strong I-scheme from $R$ to $S$.
\end{theorem}

\begin{theorem}[{\cite[Theorem 1]{aCampo_toappear_2}}] \label{theo_GschemeOnStrict}
Let $R, S \in \mfP$. Equivalent are
\begin{align*}
R & \sqsubseteq_G S; \\
\# \S(P,R) & \leq \# \S(P,S) \; \; \mytext{for all} \; P \in \mfP;  \\
\# \S(P,R) & \leq \# \S(P,S) \; \; \mytext{for all connected} \; P \in \mfP.
\end{align*}
\end{theorem}
The following proposition is not contained in \cite{aCampo_toappear_1,aCampo_toappear_2}. It shows that a mapping behaving like a (strong) Hom-scheme, G-scheme, or I-scheme on the set of the connected posets in $\mfP_r$ can be extended to a (strong) Hom-scheme, G-scheme, or I-scheme (on the total set $\mfP_r$):

\begin{proposition} \label{prop_zshP_reicht}
Let
\begin{align*}
\tau & \in \prod_{\stackrel{Q \in \mfP_r}{Q \mytext{connected}}} \A( \H(Q,R), \H(Q,S) ).
\end{align*}
For $P \in \mfP_r$ with connectivity components $K_j$, $j \in \myNk{J}$, $J \in \myN$, we define for every $\xi \in \H(P,R)$
\begin{align*}
\rho_P( \xi ) & \equiv \bigcup_{j=1}^J \tau_{K_j}( \xi \vert_{K_j} ).
\end{align*}
Then $\rho$ is a Hom-scheme, and we have:
\begin{itemize}
\item $\rho$ is strong iff $\tau_Q : \H( Q, R ) \rightarrow \H(Q,S)$ is one-to-one for every connected $Q \in \mfP_r$;
\item $\rho$ is a G-scheme iff $\mygxix{\tau_Q(\zeta)}{y} = \mygxix{\zeta}{y}$ for every connected $Q \in \mfP_r$, $\zeta \in \H(Q,R)$, and $y \in Q$;
\item $\rho$ is an I-scheme iff for every connected $Q, Q^* \in \mfP_r$, $\zeta \in \H(Q,R), \zeta^* \in \H(Q^*,R)$, $y \in Q$, $y^* \in Q^*$
\begin{equation} \label{eq_imagebased_in_Th}
\myapxix{Q}{\zeta}{y} \leq_+ \myapxix{Q^*}{\zeta^*}{y^*} \; \Rightarrow \;
\myapxix{Q}{\tau_Q(\zeta)}{y} \leq_+ \myapxix{Q^*}{\tau_{Q^*}(\zeta^*)}{y^*},
\end{equation}
where ``$<_+$'' on the left side implies ``$<_+$'' on the right side. 
\end{itemize}
\end{proposition}
\BP In all three proofs we have to show ``$\Leftarrow$'' only, because if $\tau$ fails for one of the three properties, then obviously also $\rho$ does.

It is clear that $\rho$ is a Hom-scheme. In the proofs of the first two statements, $P \in \mfP_r$ with connectivity components $K_j$, $j \in \myNk{J}$, $J \in \myN$, is fixed.

Let $\tau_Q$ be one-to-one for every connected $Q \in \mfP_r$, and let $\xi, \xi' \in \H(P,R)$ with $\xi \not= \xi'$. There exists an $x \in P$ with $\xi(x) \not= \xi'(x)$ and a $j \in \myNk{J}$ with $x \in K_j$. We conclude $\xi \vert_{K_j} \not= \xi' \vert_{K_j}$, thus $\tau( \xi \vert_{K_j} ) \not= \tau( \xi' \vert_{K_j} )$, and $\rho(\xi) \not= \rho(\xi')$ is proven.

Assume $\mygxix{\tau_Q(\zeta)}{y} = \mygxix{\zeta}{y}$ for all connected $Q \in \mfP_r$, $\zeta \in \H(Q,R)$, and $y \in Q$. Let $\xi \in \H(P,R)$ and $x \in K_j$. We have $\grxix = \gxix \subseteq K_j$, and applying \eqref{gxix_postrestr_unten} twice yields
\begin{align*}
\gxix & = \mygxix{\xi \vert_{K_j}}{x} 
= \mygxix{\tau( \xi \vert_{K_j})}{x} 
= \mygxix{\rho( \xi )\vert_{K_j}}{x} = \grxix.
\end{align*}

Now assume that \eqref{eq_imagebased_in_Th} holds for all connected posets $Q, Q^* \in \mfP_r$, $\zeta \in \H(Q,R), \zeta^* \in \H(Q^*,R)$, $y \in Q$, $y^* \in Q^*$. Let $P, P^* \in \mfP_r$, $\xi \in \H(P,R), \xi^* \in \H(P^*,R)$, $x \in P$, $x^* \in P^*$ with $\myapxix{P}{\xi}{x} \leq_+ \myapxix{P^*}{\xi^*}{x^*}$. The connectivity components of $P$ are $K_j$, $j \in \myNk{J}$, $J \in \myN$, and the connectivity components of $P^*$ are $K^*_{j^*}$, $j^* \in \myNk{J^*}$, $J^* \in \myN$. Let $j_0 \in \myNk{J}$ with $x \in K_{j_0}$ and $j^*_0 \in \myNk{J^*}$ with $x^* \in K^*_{j^*_0}$. 

Due to $ \gxix = \mygxix{\xi \vert_{K_{j_0}}}{x} \subseteq K_{j_0}$, we get
\begin{align*}
\myapxix{P}{\xi}{x}_2 & = \; \xi( \dngxix ) \; 
\stackrel{\eqref{gxix_postrestr_unten}}{=} 
\; \xi( \myodownarrow \mygxix{\xi \vert_{K_{j_0}}}{x}  )
\; = \; \xi_{K_{j_0}}( \myodownarrow \mygxix{\xi \vert_{K_{j_0}}}{x} ) \\
& = \; \myapxix{K_{j_0}}{\xi \vert_{K_{j_0}}}{x}_2
\end{align*}
and, in the same way,
$\myapxix{P}{\xi}{x}_3 = \myapxix{K_{j_0}}{\xi \vert_{K_{j_0}}}{x}_3$. Because $\xi( x ) = \xi \vert_{K_{j_0}}(x)$ is trivial, we have
$\myapxix{P}{\xi}{x} = \myapxix{K_{j_0}}{\xi \vert_{K_{j_0}}}{x}$. In the same way we see
\begin{align*}
\myapxix{P^*}{\xi^*}{x^*} & = \myapxix{K^*_{j^*_0}}{\xi^* \vert_{K^*_{j^*_0}}}{x^*}, \\
\myapxix{P}{\rho(\xi)}{x} & = \myapxix{K_{j_0}}{\tau( \xi \vert_{K_{j_0}})}{x}, \\
\myapxix{P^*}{\rho(\xi^*)}{x^*} & = \myapxix{K^*_{j^*_0}}{\tau(\xi^* \vert_{K^*_{j^*_0}})}{x^*}.
\end{align*}
Now we get
\begin{align*}
\myapxix{K_{j_0}}{\xi \vert_{K_{j_0}}}{x} & \; = \; \myapxix{P}{\xi}{x} \; \leq_+ \; \myapxix{P^*}{\xi^*}{x^*} \; = \; \myapxix{K^*_{j^*_0}}{\xi^* \vert_{K^*_{j^*_0}}}{x^*}
\end{align*}
and \eqref{eq_imagebased_in_Th} delivers
\begin{align*}
\myapxix{P}{\rho(\xi)}{x} & \; = \; \myapxix{K_{j_0}}{\tau( \xi \vert_{K_{j_0}})}{x} \; \leq_+ \; \myapxix{K^*_{j^*_0}}{\tau(\xi^* \vert_{K^*_{j^*_0}})}{x^*} \; = \; \myapxix{P^*}{\rho(\xi^*)}{x^*}
\end{align*}
with ``$<_+$'' if $\myapxix{P}{\xi}{x} <_+ \myapxix{P^*}{\xi^*}{x^*}$.

\EP

\section{Compatibility with order arithmetic} \label{sec_calcrules}

In this section, we examine to which extent the relations $\sqsubseteq, \sqsubseteq_G$, and $\sqsubseteq_I$ are compatible with order arithmetic. The simple equation \eqref{EVsys_dual} for $\E(P^d)$ lets us expect that there will be a simple relation between the different types of strong Hom-schemes from $R$ to $S$ and from $R^d$ to $S^d$. Similarly, looking at \eqref{EVsys_dirsum}, the direct sum of posets should also be easily manageable in the world of strong Hom-schemes. The formulas \eqref{EVsys_ordsum} and \eqref{EVsys_prod} for $\E(P \oplus Q)$ and $\E(P \times Q)$ refer to $\E(P)$ and $\E(Q)$ in a more involved way, and therefore, the situation will be more complicated for ordinal sums and products of posets.

\begin{proposition} \label{calc_dual}
Let $R, S \in \mfP$. Then
\begin{align*}
R \sqsubseteq S & \Leftrightarrow R^d \sqsubseteq S^d, \\
R \sqsubseteq_G S & \Leftrightarrow R^d \sqsubseteq_G S^d, \\
R \sqsubseteq_I S & \Leftrightarrow R^d \sqsubseteq_I S^d, \\
R \sqsubseteq R^d & \Leftrightarrow R \simeq R^d.
\end{align*}
\end{proposition}
\BP Let $\rho$ be a strong Hom-scheme from $R$ to $S$. Because of $\H(P,Q) = \H(P^d, Q^d)$ for every $P, Q \in \mfP$, we get a strong Hom-scheme $\rho^d$ from $R^d$ to $S^d$ by setting $\rho^d_P \equiv \rho_{P^d}$ for all $P \in \mfP$. Obviously, $\rho^d$ is a G-scheme (an I-scheme) iff $\rho$ is. The last proposition follows from the first one due to the antisymmetry of $\sqsubseteq$.

\EP

In order to avoid repetitions, we agree that the carriers of posets involved in direct or ordinal sums are always disjoint. Additionally, we agree on that the index $j$ is always an element of $\myNk{2}$.

\begin{proposition} \label{calc_dirsum}
Let $R^1, R^2, S^1, S^2 \in \mfP$. Then 
\begin{align*}
R^1 \sqsubseteq S^1 \mytext{and} R^2 \sqsubseteq S^2 & \; \Rightarrow \; R^1 + R^2 \sqsubseteq S^1 + S^2, \\
R^1 \sqsubseteq_G S^1 \mytext{and} R^2 \sqsubseteq_G S^2 & \; \Rightarrow \; R^1 + R^2 \sqsubseteq_G S^1 + S^2, \\
R^1 \sqsubseteq_I S^1 \mytext{and} R^2 \sqsubseteq_I S^2 & \; \Rightarrow \; R^1 + R^2 \sqsubseteq_I S^1 + S^2.
\end{align*}
\end{proposition}
\BP Let $\rho^j$ be a strong Hom-scheme / G-scheme / I-scheme from $R^j$ to $S^j$. For every connected $P \in \mfP$, we have $\H(P,R^1+R^2) \simeq \H(P,R^1) + \H(P,R^2)$ and $\H(P,S^1+S^2) \simeq \H(P,S^1) + \H(P,S^2)$. We define $\rho^+_P \equiv \rho^1_P \cup \rho^2_P$ for every connected $P \in \mfP$, and all statements follow with Proposition \ref{prop_zshP_reicht}. (Observe $\rho_P^+(\xi) = \rho_P^j(\xi)$, $\mygxix{\rho^+(\xi)}{x} = \mygxix{\rho^j(\xi)}{x} $, and $\myaxix{\rho^+(\xi)}{x} = \myaxix{\rho^j(\xi)}{x}$ for every $\xi \in \H(P,R^j), x \in P$, $P \in \mfP_r$ connected.)

\EP

\begin{proposition} \label{calc_ordsum}
Let $R^1, R^2, S^1, S^2 \in \mfP$. Then 
\begin{align*}
R^1 \sqsubseteq S^1 \mytext{and} R^2 \sqsubseteq S^2 & \; \Rightarrow \; R^1 \oplus R^2 \sqsubseteq S^1 \oplus S^2, \\
R^1 \sqsubseteq_G S^1 \mytext{and} R^2 \sqsubseteq_G S^2 & \; \Rightarrow \; R^1 \oplus R^2 \sqsubseteq_G S^1 \oplus S^2.
\end{align*}
\end{proposition}
\BP For posets $A, B$, the relation $A \sqsubseteq B$ means $\# \H(P,A) \leq \# \H(P,B)$ for all $P \in \mfP$, and according to Theorem \ref{theo_GschemeOnStrict}, the relation $A \sqsubseteq_G B$ is equivalent to $\# \S(P,A) \leq \# \S(P,B)$ for all $P \in \mfP$. The statements are now direct consequences of 
\begin{align*}
\H(P, A \oplus B) \simeq \sum_{U \in \U(P)} \H(X \setminus U, A) \times \H(U, B), \\
\S(P, A \oplus B) \simeq \sum_{U \in \U(P)} \S(X \setminus U, A) \times \S(U, B),
\end{align*}
for all finite posets $P, A, B$, where $\U(P)$ is the set of the upsets of $P$ and where $X \setminus U$ and $U$ are equipped with the partial order induced by $P$.

\EP

\begin{proposition} \label{calc_prod}
Let $R^1, R^2, S^1, S^2 \in \mfP$. Then 
\begin{align*}
R^1 \sqsubseteq S^1 \mytext{and} R^2 \sqsubseteq S^2 & \; \Rightarrow \; R^1 \times R^2 \sqsubseteq S^1 \times S^2, \\
R^1 \sqsubseteq_G S^1 \mytext{and} R^2 \sqsubseteq_G S^2 & \; \Rightarrow \; R^1 \times R^2 \sqsubseteq_G S^1 \times S^2.\end{align*}
\end{proposition}
\BP The first statement is a direct consequence of $\H(P, A \times B) \simeq \H(P, A) \times \H(P,B)$ for all finite posets $P, A, B$. For a constructive proof (which we need for the second statement), take strong Hom-schemes $\rho^j$ from $R^j$ to $S^j$. For every $P \in \mfP_r$ and every $\xi = ( \xi^1, \xi^2 ) \in \H( P, R^1 \times R^2 )$ (with $ \xi^j \in \H(P , R^j )$), we define $\rho^\times_P((\xi^1, \xi^2)) \equiv ( \rho^1_P(\xi^1), \rho^2_P(\xi^2))$. Then $\rho^\times$ is a strong Hom-scheme from $R^1 \times R^2$ to $S^1 \times S^2$.

Now let $\rho^1$ and $\rho^2$ be strong G-schemes. For every poset $P \in \mfP_r$, $\xi = ( \xi^1, \xi^2 ) \in \H(P,R^1 \times R^2)$, $x \in P$, we get by applying Lemma \ref{lemma_G_prod} twice
\begin{equation*}
\mygxix{\rho^\times(\xi)}{x} \; = \; \gamma_{ \mygxix{\rho^1(\xi^1)}{x} \cap \mygxix{\rho^2(\xi^2)}{x} }(x)
\; = \; \gamma_{ \mygxix{\xi^1}{x} \cap \mygxix{\xi^2}{x} }(x)
\; = \; \gxix.
\end{equation*}

\EP

\begin{proposition} \label{calc_HomSets}
Let $R, S\in \mfP$. Then 
\begin{equation} \label{HomSets_1}
R \sqsubseteq S \quad \Rightarrow \quad \H(Q,R) \sqsubseteq_G \H(Q,S) \; \mytext{for every} \; Q \in \mfP_r.
\end{equation}
Moreover,
\begin{equation} \label{HomSets_2}
\forall \; Q \in \mfP_r \; : \; \H(Q,R) \sqsubseteq \; \; \H(Q,S) \quad \Rightarrow \quad R \sqsubseteq S,
\end{equation}
and thus
\begin{align} \label{HomSets_3}
& \forall \; Q \in \mfP_r \; : \; \H(Q,R) \sqsubseteq \; \; \H(Q,S) \\
\Leftrightarrow \quad & \forall \; Q \in \mfP_r \; : \; \H(Q,R) \sqsubseteq_G \H(Q,S). \nonumber
\end{align}
Finally, 
\begin{equation} \label{HomSets_4}
\forall \; Q \in \mfP_r \; : \; \H(Q,R) \sqsubseteq_I \H(Q,S) \quad 
\Rightarrow \quad R \sqsubseteq_I S.
\end{equation}

\end{proposition}
\BP \eqref{HomSets_1}: Let $\rho$ be a strong Hom-scheme from $R$ to $S$, and let $Q \in \mfP_r$ be fixed. We define for every $P \in \mfP_r$ and every $\xi \in \H( P, \H(Q,R))$
\begin{align*}
\tau_P(\xi) & \equiv \rho_Q \circ \xi.
\end{align*}
Then $\tau_P : \H( P, \H(Q,R)) \rightarrow \H( P, \H(Q,S))$ is a well-defined mapping, and because $\rho_Q : \H(Q,R) \rightarrow \H(Q,S)$ is a one-to-one mapping, we have $\mygxix{\tau_P(\xi)}{x} = \gxix$ for all $x \in P$. $\tau$ is thus a G-scheme from $\H(Q,R)$ to $\H(Q,S)$.

Let $\xi, \zeta \in \H( P, \H(Q,R))$ with $\tau_P(\xi) = \tau_P(\zeta)$. Then $\rho_Q( \xi(x) ) = \tau_P(\xi)(x) = \tau_P(\zeta)(x) = \rho_Q( \zeta(x) )$ for all $x \in P$. Because $\rho_Q$ is one-to-one, we get $\xi(x) = \zeta(x)$ for all $x \in P$, hence $\xi = \zeta$, and $\tau$ is strong.

\eqref{HomSets_2}: In the case of $\H(Q,R) \sqsubseteq \H(Q,S) $ for all $Q \in \mfP_r$, we have for all $ P \in \mfP_r$: $ \# \H(P,R) = \# \H(A_1, \H(P,R)) \leq \# \H( A_1, \H(P,S) )$ $ = \# \H(P,S)$, thus $ R \sqsubseteq S$.

\eqref{HomSets_3}: $\H(Q,R) \sqsubseteq \H(Q,S) $ for all $Q \in \mfP_r$ yields $R \sqsubseteq S$ according to \eqref{HomSets_2}, and \eqref{HomSets_1} delivers $\H(Q,R) \sqsubseteq_G \H(Q,S) $ for all $Q \in \mfP_r$. The direction ``$\Leftarrow$'' is trivial.

\eqref{HomSets_4}: $ R \simeq \H(A_1, R) \sqsubseteq_I \H(A_1, S) \simeq S$.

\EP

\section{The ansatz for the proof of cancellation rules} \label{sec_ansatz_cancel}

In the following sections, we deal with cancellation rules for strong Hom-schemes: For $ \odot \in \{ +, \oplus, \times \}$ and ${\preceq} \in \{ \sqsubseteq, \sqsubseteq_G, \sqsubseteq_I  \}$, we have posets $Q, R, S \in \mfP$ with $ Q \odot R \preceq Q \odot S$. Does this imply $R \preceq S$? For the direct sum of posets (Section \ref{sec_cancel_dirsum}), the answer is ``Yes'' with no ifs and buts, but for ordinal sums and products (Sections \ref{sec_cancel_ordsum} and \ref{sec_cancel_product}), we need additional assumptions in most cases.

All together we prove eight cancellation rules (as already announced in the introduction, we did not succeed in proving a cancellation rule for $ Q \oplus R \sqsubseteq Q \oplus S$). In three of the proofs, the inequality between the cardinalities of the respective homomorphism sets is directly proven. The remaining five proofs use a common ansatz. Obviously, we can embed $\H(P,R)$ into $\H(P, Q \odot R)$ and $\H(P,S)$ into $\H(P, Q \odot S)$ for all $ \odot \in \{ +, \oplus, \times \}$, and \eqref{EVsys_dirsum}, \eqref{EVsys_ordsum}, and \eqref{EVsys_prod} show that also the EV-systems $\E(R)$ and $\E(S)$ can be embedded into $\E(Q \odot R )$ and $\E(Q \odot S )$, respectively, for all $ \odot \in \{ +, \oplus, \times \}$. We therefore choose a suitable embedding of $\H(P,R)$ into $\H(P, Q \odot R)$ and drive it into an embedding of $\H(P,S)$ into $\H(P, Q \odot S)$ by applying the strong Hom-scheme / G-scheme / I-scheme $\rho$ from $Q \odot R$ to $Q \odot S$. Similarly, in the case of strong I-schemes, we choose an embedding of $\E(R)$ into $\E(Q \odot R )$ and push it into an embedding of $\E(S)$ into $\E(Q \odot S)$ by means of the one-to-one homomorphism $\eps$ from $\E(Q \odot R )$ to $\E( Q \odot S)$ described in Theorem \ref{theo_strongI_equiv_eps}.

In the case of $\odot = \oplus$, the pushing process succeeds for strong I-schemes with a single application of $\eps$ on the embedding of $\E(R)$, and also for strong G-schemes, it is enough to apply $\rho$ once. In the case of $ \odot \in \{ +, \times \}$, we apply $\rho$ and $\eps$ repeatedly within iterative processes, and we show that after a finite number of iteration steps the embedding of $\H(P,R)$ arrives at an embedding of $\H(P,S)$ (the embedding of $\E(R)$ arrives at an embedding of $\E(S)$). An interesting aspect is that for both $ \odot \in \{ +, \times \}$, the iteration processes itself rely on pure set theory without order theoretical aspects. Order theory and the theory of Hom-schemes are applied on the objects the iterations are stepping through and terminating with, in showing that they provide the appropriate order theoretical properties and relations.

\section{Cancellation rules for direct sums} \label{sec_cancel_dirsum}

\begin{theorem} \label{theo_cancel_dirsum_1}
Let $Q, R, S \in \mfP$ with $Q \cap R = \emptyset$ and $Q \cap S = \emptyset$. Then
\begin{align*}
Q + R \sqsubseteq \; \; Q + S & \quad \Rightarrow \quad R \sqsubseteq S, \\
Q + R \sqsubseteq_G Q + S & \quad \Rightarrow \quad R \sqsubseteq_G S.
\end{align*}
\end{theorem}
\BP Let $P \in \mfP$ be a connected poset. For all disjoint posets $A, B \in \mfP$, we have
\begin{align*}
\H(P, A+B) & \simeq \H(P,A) + \H(P,B), \\
\S(P, A+B) & \simeq \S(P,A) + \S(P,B).
\end{align*}
$ \# \H(P, Q+R) \leq \H(P, Q+S) $ yields thus $ \# \H( P, R ) \leq \# \H( P, S )$, and with Proposition \ref{prop_zshP_reicht} we conclude $ R \sqsubseteq S$. In the same way, we get $ R \sqsubseteq_G S$ by using Theorem \ref{theo_GschemeOnStrict}.

\EP

For the proof of the cancellation rule for $\sqsubseteq_I$, we need a lemma belonging to set theory:

\begin{lemma} \label{lemma_injAbb}
Let $A, B, C$ be non-empty, disjoint sets, $A$ finite, and let the mapping $f : A \cup B \rightarrow A \cup C$ be one-to-one. For every $b \in B$, there exists an integer $n(b) \in \myN$ with
\begin{align*}
f^{n(b)}(b) & \in C, \\
f^\nu(b) & \in A \; \mytext{for} \; 1 \leq \nu \leq n(b) - 1,
\end{align*}
and the mapping
\begin{align*}
F : B & \rightarrow C, \\
b & \mapsto f^{n(b)}(b)
\end{align*}
is one-to-one.
\end{lemma}
\BP Assume $f^\nu(b) \in A$ for all $\nu \in \myN$. Because $A$ is finite, there exist $i, j \in \myN$ with $i > j$ and $f^i(b) = f^j(b)$. Because $f$ is one-to-one, we get $b = f^{i-j}(b) \in A$, a contradiction to $A \cap B = \emptyset$.

Let $x, y \in B$ with $F(x) = F(y)$, thus $f^{n(x)}(x) = f^{n(y)}(y)$. Assume $n(x) \geq n(y)$ without loss of generality. Because $f$ is one-to-one, we get $f^{n(x) - n(y)}(x) = y \in B$, and we conclude $n(x) - n(y) = 0$, thus $x = y$.

\EP

\begin{theorem} \label{theo_cancel_dirsum_2}
Let $Q, R, S \in \mfP$, with $Q \cap R = \emptyset$ and $Q \cap S = \emptyset$. Then
\begin{align*}
Q + R \sqsubseteq_I Q + S & \quad \Rightarrow \quad R \sqsubseteq_I S.
\end{align*}
\end{theorem}
\BP According to Theorem \ref{theo_strongI_equiv_eps} and \eqref{EVsys_dirsum}, there exists a one-to-one homomorphism $\eps : \E(Q) + \E(R) \rightarrow \E(Q) + \E(S)$ with $\myapexix{P}{\xi}{x} = \eps( \myapxix{P}{\xi}{x}  )$ for every $\arie$, where $\exix \; \equiv \eps( \myapxix{P}{\xi}{x} )_1.$ Due to Lemma \ref{lemma_injAbb}, there exists for every $\fa \in \E(R)$ an integer $n(\fa) \in \myN$ with $\eps^{n(\fa)}(\fa) \in \E(S)$ and
$\eps^\nu(\fa) \in \E(Q)$ for all $\nu \in \myNk{n(\fa) - 1}$; additionally, the mapping
\begin{align*}
E : \E(R) & \rightarrow \E(S), \\
\fa & \mapsto \eps^{n(\fa)}(\fa)
\end{align*}
is one-to-one. Let $\fa, \fb \in \E(R)$ with $\fa \leq_+ \fb$. Then $\fa$ and $\fb$ belong to the same connectivity component of $\E(R)$, and with $j \equiv \min\{ n(\fa), n(\fb) \}$, the points $\eps^j( \fa ) \leq_+ \eps^j( \fb )$ belong to the same connectivity component of $\E(Q) + \E(S)$. According to the definition of $j$, at least one of the points $\eps^j( \fa )$ and $\eps^j( \fb )$ belongs to $\E(S)$, thus $\eps^j( \fa ), \eps^j( \fb ) \in \E(S)$, hence $ n(\fa) =  n(\fb)$. We conclude $E( \fa ) = \eps^{n(\fa)}( \fa ) \leq_+ \eps^{n(\fa)}( \fb ) = \eps^{n(\fb)}( \fb ) = E(\fb)$, and $E$ is a one-to-one homomorphism.

Let $P \in \mfP_r, \xi \in \H(P, Q + R)$, $x \in P$. By induction, it is easily seen that $\eps^\nu( \axix ) = \myaxix{ \eta^\nu(\xi) }{x}$ for all $\nu \in \myNk{ n(\axix) }$. (Be aware that for $\nu \in \myNk{ n(\axix) - 1} $, we have $\eps^\nu( \axix ) \in \E(Q) \subseteq \E(Q) + \E(R)$, hence $\eta^\nu(\xi) \in \H(P, Q + R)$; the mapping $\eta^{\nu+1}(\xi)$ is thus well-defined.) Now we conclude $E( \axix )_1 = \eta^{n(\axix)}(\xi)(x)$ and
\begin{equation*}
E( \axix ) \; = \; \eps^{n(\axix)}(\axix) \; = \; \myaxix{ \eta^{n(\axix)}(\xi) }{x},
\end{equation*}
and $E$ fulfills the requirements of Theorem \ref{theo_strongI_equiv_eps}.

\EP

\section{Cancellation rules for ordinal sums} \label{sec_cancel_ordsum}

As already announced in the introduction, we did not succeed in proving a cancellation rule for $ Q \oplus R \sqsubseteq Q \oplus S$.

\begin{theorem} \label{theo_cancel_ordsum_G}
Let $Q, R, S \in \mfP$ with $Q \cap R = \emptyset$ and $Q \cap S = \emptyset$. Then
\begin{align*} 
Q \oplus R \sqsubseteq_G Q \oplus S & \quad \Rightarrow \quad R \sqsubseteq_G S.
\end{align*}
\end{theorem}
\BP For every $ P \in \mfP $, we define for every $\xi \in \H(P,R)$
\begin{align*}
\Gamma_{P,R}(\xi) & \; \equiv \; \mysetdescr{ \zeta \in \H(P,R) }{ \gzex = \gxix \; \mytext{for all} \; x \in P}, \\
\Gamma_{P,S}(\xi) & \; \equiv \; \mysetdescr{ \zeta \in \H(P,S) }{ \gzex = \gxix \; \mytext{for all} \; x \in P}.
\end{align*}
As proven in \cite[Lemma 4]{aCampo_toappear_2}, $R \sqsubseteq_G S$ is equivalent to
\begin{align} \label{ungl_GRxi_GSxi}
 \# \Gamma_{P,R}(\xi) & \leq \# \Gamma_{P,S}(\xi) \quad \mytext{for all} \; \xi \in \H(P,R), P \in \mfP.
\end{align}

Let $\rho$ be a strong G-scheme from $Q \oplus R $ to $Q \oplus S$, and let $P \in \mfP_r$ be fixed. For every $k \in \myN$, let $kP$ be the direct sum of $k$ disjoint isomorphic copies of $P$. We identify $\H( kP, R)$ with $\H(P, R)^k$ and $\H( kP, S)$ with $\H(P, S)^k$. Without loss of generality we can assume $Q \cap kP = \emptyset$ for every $k \in \myN$. (Otherwise, we replace $Q$ by a suitable isomorphic copy.)

Let $\xi \in \H(P,R)$ be fixed. We define for every $k \in \myN$
\begin{align*}
M_k & \equiv \mysetdescr{ \zeta \cup \chi }{ \zeta \in \S(Q,Q), \chi \in \Gamma_{P,R}(\xi)^k }
\end{align*}
$M_k$ is isomorphic to a subset of $ \H(Q,Q) \times \H( kP, R )$ which is again isomorphic to a subset of $\H(Q \oplus kP, Q \oplus R )$. $\rho_k( \theta ) \equiv \rho_{Q \oplus kP}(\theta)$ is thus well-defined for every $\theta \in M_k$.

Let $\theta = \zeta \cup \chi \in M_k$ with $\zeta \in \S(Q,Q)$, $\chi \in \Gamma_{P,R}( \xi )^k$. We have $\theta(Q) \subseteq Q$ and $\theta( kP ) \subseteq R$ with $Q \cap R = \emptyset$, which yields $\mygxix{\theta}{y} = \mygxix{\zeta}{y}$ for $y \in Q$. Therefore, for $y \in Q$, we get $\mygxix{\rho_k(\theta)}{y} = \mygxix{\theta}{y} = \{ y \}$, because $\zeta$ is strict. We conclude $ \rho_k(\theta) \vert_Q \in \S(Q, Q \oplus S)$.

Now we show $\rho_k(\theta)( k P ) \subseteq S$. Let $y_1 < \ldots < y_L$ be a chain of maximal length in $Q$. Due to the strictness of $\rho_k(\theta) \vert_Q$, we have $\rho_k(\theta)(y_1) < \ldots < \rho_k(\theta)(y_L)$. Assume $\rho_k(\theta)(x) \in Q$ for an $x \in kP$. Due to $y_\ell < x$ in $Q \oplus kP$ for all $\ell \in \myNk{L}$, we have $\rho_k(\theta)(y_\ell) \leq \rho_k(\theta)(x)$, and, in particular, $\rho_k(\theta)(y_\ell) \in Q$ for all $\ell \in \myNk{L}$. We conclude $\rho_k(\theta)(y_L) = \rho_k(\theta)(x)$, because in the case of ``$<$'' we get a chain of length $L+1$ in $Q$. But in the case $\rho_k(\theta)(y_L) = \rho_k(\theta)(x)$, we have $x \in \mygxix{\rho_k(\theta)}{y_L} = \{ y_L \}$ in contradiction to $Q \cap kP = \emptyset$. Thus, $\rho_k(\theta)( k P ) \subseteq S$. In particular, $\rho_k(\theta) \vert_{kP}^S$ is well-defined.

Let $x \in kP$. Due to $\theta(Q) \subseteq Q$ and $\theta( kP ) \subseteq R$ with $Q \cap R = \emptyset$, we have $ \mygxix{\theta}{x} \subseteq kP $, hence $ \mygxix{\rho_k(\theta)}{x} \subseteq kP $, too. Now we get
\begin{align*}
\mygxix{\rho_k(\theta) \vert_{kP}^S}{x} & 
\stackrel{(\ref{gxix_postrestr_oben},\ref{gxix_postrestr_unten})}{=}
\mygxix{\rho_k(\theta)}{x} = \mygxix{\theta}{x} 
\stackrel{\eqref{gxix_postrestr_unten}}{=}
\mygxix{\theta \vert_{kP}}{x}
= \mygxix{\chi}{x}
= \gxix,
\end{align*}
thus $\rho_k(\theta) \vert_{kP}^S \in \Gamma_{P,S}( \xi )^k$. All together, we can write $\rho_k(\theta) = \zeta' \cup \chi'$ with $\zeta' \in \S( Q, Q \oplus S)$ and $\chi' \in \Gamma_{P,S}( \xi )^k$. Accordingly, we can regard $\rho_k( M_k )$ as a subset of $\S( Q, Q \oplus S) \times \Gamma_{P,S}( \xi )^k$. Now we get
\begin{align*}
\# \S(Q,Q) \cdot ( \# \Gamma_{P,R}( \xi ) )^k
& = \# ( \S(Q, Q) \times \Gamma_{P,R}( \xi )^k ) \\
& = \# M_k \; = \; \# \rho_k( M_k ) \\
& \leq \# ( \S(Q, Q \oplus S) \times \Gamma_{P,S}( \xi )^k ) \\
& = \# \S(Q,Q \oplus S) \cdot ( \# \Gamma_{P,S}( \xi ) )^k.
\end{align*}
In the case $\# \Gamma_{P,R}( \xi ) = 1$, \eqref{ungl_GRxi_GSxi} is now trivial, and in the case $\# \Gamma_{P,R}( \xi ) > 1$, \eqref{ungl_GRxi_GSxi} holds because $k \in \myN$ is arbitrary. 

\EP

The rest of the section is dedicated to the proof of a cancellation rule for I-schemes. In what follows, $Q$, $R$, and $S$ are fixed finite non-empty posets with $Q \oplus R \sqsubseteq_I Q \oplus S$, and the one-to-one homomorphism $\eps : \E( Q \oplus R ) \rightarrow \E( Q \oplus  S )$ has the property described in Theorem \ref{theo_strongI_equiv_eps}. $X$ is the carrier of $R$ and $Y$ is the carrier of $Q$. Furthermore:
\begin{align*}
\E(Q)^* & \equiv \mysetdescr{ (x, D, U \cup X) }{(x,D,U) \in \E(Q) }, \\
\E(R)^* & \equiv \mysetdescr{ (x, D \cup Y, U ) }{(x,D,U) \in \E(R) }, \\
\E(S)^* & \equiv \mysetdescr{ (x, D \cup Y, U ) }{(x,D,U) \in \E(S) }.
\end{align*}
According to \eqref{EVsys_ordsum}, $\E(Q)^*, \E(R)^* \subseteq \E( Q \oplus R )$ with $\fa <_+ \fb $ for every $\fa \in \E(Q)^*, \fb \in \E(R)^*$.

\begin{lemma} \label{lemma_epsER_epsEQ}
We have
\begin{equation} \label{eq_epsER}
\eps( \E(R)^* )_1 \; \subseteq \; S,
%\mysetdescr{ \fa \in \E( Q \oplus S ) }{ \fa_1 \in S },
\end{equation}
and in the case $Y \subseteq \eps( \E(Q)^* )_1$ we have
\begin{equation} \label{eq_epsERES}
\eps( \E(R)^* ) \; \subseteq \; \E(S)^*.
\end{equation}
\end{lemma}
\BP Obviously, $h_{\E(Q)^*} = h_{\E(Q)} = h_Q$. Let $\fa^1 <_+ \ldots <_+ \fa^{h_Q}$ be a chain of maximal length in $\E(Q)^*$ and let $\fb \in \E(R)^*$, hence $\fa^{h_Q} <_+ \fb$. Then $\eps(\fa^1) <_+ \ldots <_+ \eps(\fa^{h_Q}) <_+ \eps(\fb)$ is a chain in $\E(Q \oplus S)$, and $\eps(\fa^1)_1 < \ldots < \eps(\fa^{h_Q})_1 < \eps(\fb)_1$ is a chain in $Q \oplus S$. We conclude $\eps(\fb)_1 \in S$.

Now assume $Y \subseteq \eps( \E(Q)^* )_1$, and let $\fb \in \E(R)^*$ be fixed. According to \eqref{eq_epsER}, $\eps( \fb )_1 \in S$. Furthermore, due to $\fa <_+ \fb$ for all $\fa \in \E(Q)^*$, we have $\eps( \fa )_1 \in \eps( \fb )_2$ for every $\fa \in \E(Q)^*$. According to our assumption, this means $Y \subseteq \eps( \fb )_2$, hence $\eps( \fb ) \in \E( S )^*$.

\EP

We now define two mappings:
\begin{align*}
r : \E(R) \quad & \rightarrow \E(Q \oplus R), \\
( x, D, U ) & \mapsto ( x, D \cup Y, U ); \\
s : \E( Q \oplus S ) & \rightarrow \E(S), \\
( x, D, U ) & \mapsto ( x, D \setminus Y, U ).
\end{align*}

\begin{corollary} \label{corr_sepsr_is_hom}
Assume $Y \subseteq \eps( \E(Q)^* )_1$. Then
\begin{equation*}
s \circ \eps \circ r : \E(R) \rightarrow \E(S)
\end{equation*}
is a one-to-one homomorphism.
\end{corollary}
\BP Obviously, $r$ is an embedding with $r( \E(R) ) = \E(R)^*$. For every $\fa \in \E(R)$, equation \eqref{eq_epsERES} yields $ \eps( r( \fa ) ) \in \E(S)^* $, and because $s$ is one-to-one on $\E(S)^*$, the total mapping $s \circ \eps \circ r$ is one-to-one.

Let $\fa, \fa' \in \E(R)$ with $ \fa <_+ \fa'$. Because $r$ is an embedding and $\eps$ a strict homomorphism, we have $\eps(r(\fa)) <_+ \eps(r(\fa'))$, hence $\eps(r(\fa))_1 \in \eps(r(\fa'))_2$ and $\eps(r(\fa'))_1 \in \eps(r(\fa))_3$. Due to $\eps(r(\fa)), \eps(r(\fa')) \in \E(S)^*$, there exist $ (x , D, U )$, $(x', D', U' ) \in \E(S)$ with $( x, D \cup Y, U ) = \eps( r( \fa ) )$ and $( x', D' \cup Y, U' ) = \eps( r( \fa' ) )$. We conclude $ x \in D' $ and $ x' \in U $, thus $s(\eps(r( \fa ))) = ( x, D, U ) <_+ ( x', D', U' ) = s(\eps(r( \fa' )))$. Hence, $s \circ \eps \circ r$ is a homomorphims.

\EP

\begin{lemma} \label{lemma_epsr_23}
Assume $Y \subseteq \eps( \E(Q)^* )_1$. Then
\begin{align} \label{eq_epsraxix_2}
\eps( r ( \axix ) )_2 & = 
\mysetdescr{ \eps( r ( \axiy ) )_1 }{ y \in \dngxix } \; \cup \; Y ,\\ \label{eq_epsraxix_3}
\eps( r ( \axix ) )_3 & = 
\mysetdescr{ \eps( r ( \axiy ) )_1 }{ y \in \upgxix },
\end{align}
for all $\arie$. In \eqref{eq_epsraxix_2}, the two sets on the right are disjoint.
\end{lemma}
\BP Let $\arie$ be fixed. Due to $r ( \axix ) = \myapxix{Q \oplus P}{ id_Q \cup \xi}{x}$ we get with $\eta$ defined as in Theorem \ref{theo_strongI_equiv_eps}
\begin{align*}
\eps( r ( \axix ) ) & = \eps( \myapxix{Q \oplus P}{ id_Q \cup \xi}{x} ) \; = \; \myapxix{Q \oplus P}{ \eta( id_Q \cup \xi )}{x}.
\end{align*}
The mapping $\eta$ is an I-scheme and thus a G-scheme. Therefore, $\mygxix{\eta(\zeta)}{z} = \mygxix{\zeta}{z}$ for all $\zeta \in \H(Q \oplus P, Q \oplus R)$, $z \in Q \oplus P$, and thus 
\begin{align*}
\eps( r ( \axix ) )_2 & = \myapxix{Q \oplus P}{ \eta( id_Q \cup \xi )}{x}_2 \\
& = \eta( id_Q \cup \xi )\left( \mydngxix{id_Q \cup \xi}{x} \right) \\
& = \mysetdescr{ \eta( id_Q \cup \xi )(y) }{ y \in \mydngxix{id_Q \cup \xi}{x} } \\
& = \mysetdescr{ \eps( \myapxix{Q \oplus P}{ id_Q \cup \xi}{y} )_1 }{ y \in \mydngxix{id_Q \cup \xi}{x} }.
\end{align*}
Because of $x \in R$, we have ${\myodownarrow_{Q \oplus P}} \mygxix{id_Q \cup \xi}{x} = Y \cup {\myodownarrow_{P}} \gxix$. With
\begin{align*}
N_Q & \equiv \mysetdescr{ \eps( \myapxix{Q \oplus P}{ id_Q \cup \xi}{y} )_1 }{ y \in Y }, \\
N_P & \equiv \mysetdescr{ \eps( \myapxix{Q \oplus P}{ id_Q \cup \xi}{y} )_1 }{ y \in {\myodownarrow_{P}} \gxix }
\end{align*}
we have thus $\eps( r( \axix ) )_2 = N_Q \cup N_P$.

Let $y \in {\myodownarrow_{P}} \gxix$, hence $y \in P$. Then we have $Y \subseteq  \myapxix{Q \oplus P}{ id_Q \cup \xi}{y}_2$, thus $ \myapxix{Q \oplus P}{ id_Q \cup \xi}{y} \in \E(R)^*$. Equation \eqref{eq_epsER} delivers $\eps( \myapxix{Q \oplus P}{ id_Q \cup \xi}{y})_1 \in S$, and we conclude $N_P \subseteq S$. In particular, $N_P \cap Y = \emptyset$.

Due to \eqref{eq_epsERES}, we have $ Y \subseteq \eps( r( \axix ) ) )_2$. We conclude $Y \subseteq N_Q$, and with $ \# N_Q \leq \# Y$ we get $Y = N_Q$. Now \eqref{eq_epsraxix_2} is shown.

In the same way we see
\begin{align*}
\eps( r ( \axix ) )_3 & = \mysetdescr{ \eps( \myapxix{Q \oplus P}{ id_Q \cup \xi}{y} )_1 }{ y \in \myupgxix{id_Q \cup \xi}{x} }
\end{align*}
and \eqref{eq_epsraxix_3} follows because of ${{\myouparrow}_{Q \oplus P}} \mygxix{id_Q \cup \xi}{x} = {{\myouparrow}_{P}} \gxix$.

\EP

\begin{theorem} \label{theo_cancel_ordsumI}
Assume $Y \subseteq \eps( \E(Q)^* )_1$. Then $R \sqsubseteq_I S$. In particular, the mapping 
\begin{align*}
T \equiv s \circ \eps \circ r : \E(R) \rightarrow \E(S)
\end{align*}
is a one-to-one homomorphism which fulfills the requirement of Theorem \ref{theo_strongI_equiv_eps}: for all $\arie$
\begin{align*}
\myaxix{\tau_P(\xi)}{x} & = T( \myapxix{P}{\xi}{x} ), \\
\mytext{where} \quad
\tau_P(\xi)(x) & \equiv T( \myapxix{P}{\xi}{x} )_1.
\end{align*}
\end{theorem}
\BP In Corollary \ref{corr_sepsr_is_hom}, we have seen that $T$ is a one-to-one homomorphism. 
$\myaxix{\tau(\xi)}{x}_1 = T( \axix )_1$ is trivial. Due to Lemma \ref{lemma_exix},  $\mygxix{\tau(\xi)}{x} = \gxix$, and with $\tau(\xi)(x) = T( \axix )_1 = \eps( r( \axix ) )_1$ we get
\begin{align*}
\myaxix{\tau(\xi)}{x}_2 & = \tau(\xi)( \dngxix ) = \mysetdescr{ \eps( r ( \axiy ) )_1 }{ y \in \dngxix } \\
& \stackrel{\eqref{eq_epsraxix_2}}{=} \eps( r ( \axix ) )_2 \setminus Y = s( \eps( r( \axix ) ) )_2 = T( \axix )_2.
\end{align*}
$ \myaxix{\tau(\xi)}{x}_3 = T( \axix )_3$ is shown in the same way by using \eqref{eq_epsraxix_3}.

\EP

The dual results are

\begin{theorem} \label{theo_cancel_ordsum_duals}
Let $Q, R, S \in \mfP$ with $Q \cap R = \emptyset$ and $Q \cap S = \emptyset$. Then
\begin{equation*}
R \oplus Q \sqsubseteq_G S \oplus Q \quad \Rightarrow \quad R \sqsubseteq_G S.
\end{equation*}
In the case $R \oplus Q \sqsubseteq_I S \oplus Q$, let $\eps : \E( R \oplus Q ) \rightarrow \E( S \oplus Q )$ be the one-to-one homomorphism from Theorem \ref{theo_strongI_equiv_eps}, and let
\begin{equation*}
\E(Q)_* \equiv \mysetdescr{ (x, D \cup X, U) }{(x,D,U) \in \E(Q) },
\end{equation*}
where $X$ is the carrier of $R$. With $Y$ being the carrier of $Q$, we have
\begin{equation*}
Y \subseteq \eps( \E(Q)_* )_1 \quad \Rightarrow \quad R \sqsubseteq_I S.
\end{equation*}
\end{theorem}

\section{Cancellation rules for products} \label{sec_cancel_product}

In this section, $Q, R, S \in \mfP$ are fixed posets, and $\rho$ is a strong Hom-scheme from $Q \times R$ to $Q \times S$. The first cancellation rule is easily proven: For all $P \in \mfP_r$ we have
\begin{align*}
\# \H(P,Q) \cdot \# \H(P,R) & = \# \H(P, Q \times R) \; \leq \; \# \H(P, Q \times S) \\
& = \# \H(P,Q) \cdot \# \H(P,S),
\end{align*}
thus $R \sqsubseteq S$.

In the proofs of the cancellation rules for strong G-schemes and I-schemes, constant homomorphisms will play an important role:
\begin{definition} \label{def_CP_cPq}
For all $P \in \mfP$ and all $q \in Q$, let $c_{P,q} \in \H( P, Q )$ be the homomorphism mapping all $x \in P$ to $q$.
\begin{equation*}
\C(P) \; \equiv \; \mysetdescr{ c_{P,q} \in \H(P,Q) }{ q \in Q }
\end{equation*}
is the set of the constant homomorphims from $P$ to $Q$.
\end{definition}

\begin{corollary} \label{corr_cqxi}
Let $P \in \mfP$, $q \in Q$, and $\xi \in \H(P,R)$. Then for all $x \in P$
\begin{align} \label{eq_gxix_cqxi}
\mygxix{(c_{P,q},\xi)}{x} & = \gxix, \\
\myaxix{(c_{P,q},\xi)}{x} & = \label{eq_axix_cqxi}
\big( ( q, \xi(x) ),
\{ q \} \times \axix_2,
\{ q \} \times \axix_3 \big).
\end{align}
\end{corollary}
\BP $\mygxix{c_{P,q}}{x}$ is for every $x \in P$ the connectivity component of $P$ containing $x$, thus $\gxix \subseteq \mygxix{c_{P,q}}{x}$ because $\gxix$ is connected with $x \in \gxix$. Hence,
\begin{align*}
\mygxix{ ( c_{P,q}, \xi) }{x} & \stackrel{\eqref{gxix_produkt}}{=} 
\gamma_{\mygxix{c_{P,q}}{x} \cap \gxix}(x)
= \gamma_{\gxix}(x) \stackrel{\eqref{gxix_gammagxix}}{=} \gxix.
\end{align*}
$ \myaxix{(c_{P,q},\xi)}{x}_1 = (q, \xi(x) )$ is trivial. Furthermore,
\begin{align*}
\myaxix{(c_{P,q},\xi)}{x}_2 & = 
( c_{P,q}, \xi)\Big( {\myodownarrow} \mygxix{ ( c_{P,q}, \xi) }{x} \Big)
\stackrel{\eqref{eq_gxix_cqxi}}{=}
( c_{P,q}, \xi)\Big( {\myodownarrow} \gxix \Big) \\
& = 
\{ q \} \times \xi\Big( {\myodownarrow} \gxix \Big) = 
\{ q \} \times \axix_2
\end{align*}
and similarly $\myaxix{(c_{P,q},\xi)}{x}_3 = \{ q \} \times \axix_3$.

\EP

In this section, the key is the following set-theoretical lemma:

\begin{lemma} \label{lemma_produktzyklus}
Let $A, B, C$ be non-empty sets, $A$ finite, and let $f : A \times B \rightarrow C$ and $g : C \rightarrow A$ be mappings with
\begin{equation} \label{bed_komp_1_allg}
\forall \; a, a' \in A, b \in B \; : \quad a \not= a' \; \Rightarrow \; g(f(a,b)) \not= g(f(a',b)).
\end{equation}
For a fixed $a_0 \in A$, we define recursively for every $b \in B$ the sequence
\begin{align*}
\chi^0 (b) & \equiv a_0, \\
\forall \; i \in \myN \; : \; \chi^i (b) & \equiv g( f ( \chi^{i-1}(b), b)).
\end{align*}
Then, for every $b \in B$, there exists an integer $n(b) \in \myN$ with
\begin{align*}
\chi^{n(b)} (b) & = a_0.
\end{align*}
\end{lemma}
\BP Let $b \in B$. We have $\chi^i(b) \in A$ for every $i \in \myN_0$, and because $A$ is finite, there exist $i, j \in \myN_0$ with $i > j$ and $\chi^i(b) = \chi^j(b)$. In the case $j = 0$, we have $\chi^i(b) = a_0$. For $j > 0$ we get
\begin{align*}
g(f( \chi^{i-1}(b), b)) & = \; \chi^i(b) \spc{=} \chi^j(b) \spc{=} g(f( \chi^{j-1}(b), b )),
\end{align*}
and with \eqref{bed_komp_1_allg} we conclude $\chi^{i-1}(b) = \chi^{j-1}(b)$. Stepping backwards in this way, we finally reach $\chi^{i-j}(b) = \chi^0(b) = a_0$.

\EP

Now let $\rho$ be a G-scheme from $Q \times R$ to $Q \times S$. In order to be able to show $R \sqsubseteq_G S$, we will assume that $\rho$ fulfills \eqref{bed_komp_1_allg}:
\begin{align} \label{bed_komp_1}
\forall P \in \mfP,  \xi, \zeta \in \H(P,Q), \theta \in \H(P,R) & :
\\
\xi \not= \zeta \; & \Rightarrow \; \rho( \xi, \theta )_1 \not= \rho( \zeta, \theta )_1, \nonumber
\end{align}
thus, (with the symbols used in Lemma \ref{lemma_produktzyklus}): $A = \H(P,Q)$, $B = \H(P,R)$, $C = \H(P, Q \times S)$, $f = \rho_P$, $g = \pi_1$. The counterpart to the sequences $\chi^i(b)$ is defined as follows:
\begin{definition} \label{def_phi}
Let $q \in Q$ be fixed. For every $\xi \in \H(P, R )$, we define recursively 
\begin{align*}
\phi_P^0(\xi) & \equiv c_{P,q} \\
\forall \; i \in \myN \; : \; 
\phi_P^{i}(\xi) & \equiv \rho( \phi_P^{i-1}(\xi), \xi )_1.
\end{align*}
As complementary notation, we define
\begin{align*}
\forall \; i \in \myN_0 \; : \; 
\Phi_P^{i}(\xi) & \equiv ( \phi_P^i(\xi), \xi ) = 
\begin{cases}
( c_{P,q}, \xi ), & \mytext{if} i = 0; \\
\left( \rho( \Phi_P^{i-1}(\xi) )_1, \xi \right), & \mytext{if} i \geq 1.
\end{cases}
\end{align*}
\end{definition}

\begin{lemma} \label{lemma1_phi_zu_QxR}
For every $\arie$, and every $i \in \myN_0$
\begin{align} \label{grhophix_eq_gxix_1}
\mygxix{ \Phi_P^i(\xi) }{x} & = \gxix.
\end{align}
\end{lemma}
\BP Let $\arie$ be selected. In order to unburden the notation, we skip the index $P$. $\mygxix{ \Phi^0(\xi) }{x} = \gxix$ holds due to Corollary \ref{corr_cqxi}.

Assume that \eqref{grhophix_eq_gxix_1} holds for $i \in \myN_0$. Then $\gxix = \mygxix{ \Phi^i(\xi) }{x}$ $= \mygxix{ \rho( \Phi^i(\xi)) }{x}$ $\subseteq \mygxix{\rho( \Phi^{i}(\xi) )_1}{x}$ according to \eqref{gxix_produkt}, thus
\begin{align*}
\mygxix{ \Phi^{i+1}(\xi) }{x} & = 
\mygxix{ ( \rho( \Phi^{i}(\xi) )_1, \xi) }{x}
\stackrel{\eqref{gxix_produkt}}{=} 
\gamma_{\mygxix{\rho( \Phi^{i}(\xi))_1}{x} \; \cap \; \gxix}(x) \\
& = \gamma_{\gxix}(x) \stackrel{\eqref{gxix_gammagxix}}{=} \gxix.
\end{align*}

\EP

\begin{theorem} \label{theo_cancel_prod_G}
Let $\rho$ fulfill \eqref{bed_komp_1}. According to Lemma \ref{lemma_produktzyklus}, there exists for every $\xi \in \H(P,R)$ an integer $n_P(\xi) \in \myN$ with
\begin{align*}
\phi_P^{n_P(\xi)}(\xi) & = c_{P,q}.
\end{align*}
Then the mapping
\begin{align*}
\tau & \in \prod_{P \in \mfP_r} \A( \H(P,R), \H(P,S) ) \\
\mytext{defined by} \quad
\tau_P(\xi) & \equiv \rho( \Phi_P^{n_P(\xi)-1}(\xi) )_2
\end{align*}
for all $P \in \mfP_r, \xi \in \H(P,R)$, is a strong G-scheme from $R$ to $S$.
\end{theorem}
\BP Let $\arie$. Skipping the index $P$, we have 
$\rho( \Phi^{n(\xi)-1}(\xi))_1 = \phi^{n(\xi)}(\xi) = c_{P,q}$, thus
\begin{align*}
\mygxix{ \tau( \xi ) }{x} & 
\stackrel{\eqref{eq_gxix_cqxi}}{=}
\mygxix{ \rho( \Phi^{n(\xi)-1}(\xi)) }{x} = 
\mygxix{ \Phi^{n(\xi)-1}(\xi) }{x} 
\stackrel{\eqref{grhophix_eq_gxix_1}}{=} \gxix,
\end{align*}
and $\tau$ is a G-scheme.

Let $\xi, \zeta \in \H(P,R)$ with $\tau(\xi) = \tau(\zeta)$. Then
\begin{align*}
\rho( \Phi^{n(\xi)-1}(\xi)_1, \xi ) & = 
( c_{P,q}, \tau(\xi) )
= ( c_{P,q}, \tau(\zeta) ) 
= \rho( \Phi^{n(\zeta)-1}(\zeta)_1, \zeta ),
\end{align*}
and because $\rho$ is one-to-one we get $ ( \Phi^{n(\xi)-1}(\xi)_1, \xi ) = ( \Phi^{n(\zeta)-1}(\zeta)_1, \zeta )$, thus $\xi = \zeta$, and $\tau_P$ is one-to-one.

\EP

Also in this section, the proof of the cancellation rule for I-schemes will be the more complicated one. Assume $ Q \times R \sqsubseteq_I Q \times S $, and let $\rho$ be a strong I-scheme from $R$ to $S$. For the proof of $R \sqsubseteq_I S$, we need an additional assumption about $\rho$:
\begin{align} \label{assumption_Ischeme}
\forall \; \xi \in \H(P, Q \times R) & \; : \;
\xi_1 \in \C(P) \; \Rightarrow \; \rho_P(\xi)_1 \in \C(P).
\end{align}

\begin{lemma} \label{lemma_alpha_Phi}
Let $\rho$ fulfill \eqref{assumption_Ischeme}. Then for all $P, P' \in \mfP_r, \xi \in \H(P, Q \times R)$, $\zeta \in \H(P', Q \times R)$, $x \in P$, $y \in P'$
\begin{align} \label{aPhiixi_le_aPhiizeta}
\myapxix{P}{\xi}{x} \leq_+ \myapxix{P'}{\zeta}{y} & \; \Rightarrow \;
\forall \; i \in \myN_0 \; : \; \myapxix{P}{\Phi_P^i(\xi)}{x} \leq_+ \myapxix{P'}{\Phi_{P'}^i(\zeta)}{y},
\end{align}
where ``$<_+ / = $'' on the left side implies ``$<_+ / =$'' on the right side. Additionally,
\begin{align} \label{Phiixi_eq_Phiizeta}
\myapxix{P}{\xi}{x} \leq_+ \myapxix{P'}{\zeta}{y} & \; \Rightarrow \;
\forall \; i \in \myN_0 \; : \; \Phi_P^i(\xi)_1 = \Phi_{P'}^i(\zeta)_1
\end{align}
for every $P, P' \in \mfP_r$, $\xi \in \H(P, Q \times R)$, $\zeta \in \H(P', Q \times R)$, $x \in P$, and $y \in P'$.
\end{lemma}
\BP Also in this proof, we skip the indices $P$ and $P'$. For $\axix = \myaxix{\zeta}{y}$, \eqref{eq_axix_cqxi} delivers $\myaxix{\Phi^0(\xi)}{x} = \myaxix{\Phi^0(\zeta)}{y}$. Assume $\axix <_+ \myaxix{\zeta}{y}$. Then $\xi(x) \in \myaxix{\zeta}{y}_2$ and $\zeta(y) \in \myaxix{\xi}{x}_3$, thus
\begin{align*}
\Phi^0(\xi)(x) & = ( q, \xi(x) ) \in \{ q \} \times \myaxix{\zeta}{y}_2 
\stackrel{\eqref{eq_axix_cqxi}}{=} \myaxix{\Phi^0(\zeta)}{y}_2
\end{align*}
and similarly $ \Phi^0(\zeta)(y) \in \myaxix{\Phi_P^0(\xi)}{x}_3$. Therefore, $\myaxix{\Phi^0(\xi)}{x} <_+ \myaxix{\Phi^0(\zeta)}{y}$, and \eqref{aPhiixi_le_aPhiizeta} holds for $i = 0$. \eqref{Phiixi_eq_Phiizeta} is trivial for $i = 0$. 

Assume that $\myaxix{\Phi^i(\xi)}{x} \leq_+ \myaxix{\Phi^i(\zeta)}{y}$ holds for $ i \in \myN_0 $. Because $\rho$ is an I-scheme, we have $\myarxix{\Phi^i(\xi)}{x}$ $\leq_+$ $\myarxix{\Phi^i(\zeta)}{y}$. By applying \eqref{assumption_Ischeme} $(i+1)$ times we see that $ \rho(\Phi^i(\xi))_1 $, and $ \rho(\Phi^i(\zeta))_1 $ are both elements of $\C(P)$. There exist thus $r, s \in Q$ with $ \rho(\Phi^i(\xi))_1 = c_{P,r} $ and $ \rho(\Phi^i(\zeta))_1 = c_{P',s} $.

In the case $\axix = \myaxix{\zeta}{y}$, we have $\myaxix{\Phi^i(\xi)}{x} = \myaxix{\Phi^i(\zeta)}{y}$, and \eqref{eq_imagebased_rhowert} yields $\myarxix{\Phi^i(\xi)}{x}$ $=$ $\myarxix{\Phi^i(\zeta)}{y}$. Now
\begin{align*}
( r, \rho(\Phi^i(\xi))_2(x) ) & = \rho(\Phi^i(\xi))(x) = \rho(\Phi^i(\zeta))(y) = ( s, \rho(\Phi^i(\zeta))_2(y) ),
\end{align*}
thus $r = s$. And in the case $\axix <_+ \myaxix{\zeta}{y}$, we have $\myaxix{\Phi^i(\xi)}{x} <_+ \myaxix{\Phi^i(\zeta)}{y}$, hence $\myarxix{\Phi^i(\xi)}{x}$ $<_+$ $\myarxix{\Phi^i(\zeta)}{y}$. Therefore,
\begin{align*}
( r, \rho(\Phi^i(\xi))_2(x) ) & =  \rho(\Phi^i(\xi))(x) \in
\myarxix{\Phi^i(\zeta)}{y}_2 \stackrel{\eqref{eq_axix_cqxi}}{=}
 \{ s \} \times \myaxix{\rho(\Phi^i(\zeta))_2}{y}_2,
\end{align*}
hence $r = s$ also in this case. \eqref{Phiixi_eq_Phiizeta} holds thus for $i+1$, and we get with $\xi(x) \in \myaxix{\zeta}{y}_2$
\begin{align*}
\Phi^{i+1}(\xi)(x) & = ( r, \xi(x) ) \in \{ r \} \times \zeta \left( \myodownarrow \mygxix{\zeta}{y} \right) = 
\{ s \} \times \zeta \left( \myodownarrow \mygxix{\zeta}{y} \right) \\
& = ( c_{P,s}, \zeta ) \left( \myodownarrow \mygxix{\zeta}{y} \right)
\stackrel{\eqref{grhophix_eq_gxix_1}}{=}
\myaxix{\Phi^{i+1}(\zeta)}{y}_2,
\end{align*}
and in the same way we prove $ \Phi^{i+1}(\zeta)(y) \in \myaxix{\Phi^{i+1}(\xi)}{x}_3$. Thus, $\myaxix{\Phi^{i+1}(\xi)}{x}$ $<_+$ $\myaxix{\Phi^{i+1}(\zeta)}{y}$.

\EP

\eqref{Phiixi_eq_Phiizeta} shows that assumption \eqref{assumption_Ischeme} is restrictive. For all $i \in \myN_0$, it enforces $\phi_P^i(\xi) = \phi_P^i(\zeta)$ for all $\xi, \zeta \in \H(P, Q \times R)$ for which $\myapxix{P}{\xi}{x}$ and $\myapxix{P}{\zeta}{y}$ belong to the same connectivity component of $\E(Q \times R)$ for some $x, y \in P$.

\begin{theorem} \label{theo_cancel_prod_I}
Let $\rho$ fulfill \eqref{bed_komp_1} and \eqref{assumption_Ischeme}.
Then $R \sqsubseteq_I S$, and the mapping $\tau$ defined in Theorem \ref{theo_cancel_prod_G} is a strong I-scheme from $R$ to $S$.
\end{theorem}
\BP Let $P \in \mfP_r, \xi \in \H(P,Q \times R)$, and let $n_P(\xi)$ be defined as in Theorem \ref{theo_cancel_prod_G}. Due to $\phi^{n_P(\xi)}( \xi ) = c_{P,q}$, we have $ \Phi^{n_P(\xi)}( \xi ) = ( c_{P,q}, \xi ) = \Phi^0( \xi )$. We conclude $ \Phi^{k \cdot n_P(\xi) - 1}( \xi ) = \Phi^{n_P(\xi) - 1}( \xi )$ for all $k \in \myN$.

Let $P, P' \in \mfP_r, \xi \in \H(P, Q \times R)$, $\zeta \in \H(P', Q \times R)$, $x \in P$, $y \in P'$ with $\myaxix{\xi}{x} \leq_+ \myaxix{\zeta}{y}$. (Now we skip again the indices $P$ and $P'$.) With $N \equiv \lcm \{n(\xi), n(\zeta) \}$, Lemma \ref{lemma_alpha_Phi} yields $\myaxix{\Phi^{N-1}(\xi)}{x} \leq_+ \myaxix{\Phi^{N-1}(\zeta)}{y}$ with ``$<_+$'' if $\myaxix{\xi}{x} <_+ \myaxix{\zeta}{y}$ and with ``$=$'' if $\myaxix{\xi}{x} = \myaxix{\zeta}{y}$. Because $\rho$ is an I-scheme, we get 
\begin{align*}
\myarxix{\Phi^{n(\xi)-1}(\xi)}{x}
& =  \myarxix{\Phi^{N-1}(\xi)}{x}
\leq_+ \myarxix{\Phi^{N-1}(\zeta)}{y}
= \myarxix{\Phi^{n(\zeta)-1}(\zeta)}{y},
\end{align*}
again with ``$<_+$'' if $\myaxix{\xi}{x} <_+ \myaxix{\zeta}{y}$, and with ``$=$'' if $\myaxix{\xi}{x} = \myaxix{\zeta}{y}$, according to \eqref{eq_imagebased_rhowert}.

Assume $\myaxix{\xi}{x} = \myaxix{\zeta}{y}$. Then, with Corollary \ref{corr_cqxi}
\begin{align*}
& \; \big( ( q, \myaxix{\tau(\xi)}{x}_1 ),
\{ q \} \times \myaxix{\tau(\xi)}{x}_2,
\{ q \} \times \myaxix{\tau(\xi)}{x}_3 \big) \\
= & \; \myarxix{\Phi^{n(\xi)-1}(\xi)}{x}
= \myarxix{\Phi^{n(\xi)-1}(\zeta)}{y} \\
= &  \; \big( ( q, \myaxix{\tau(\zeta)}{y}_1 ),
\{ q \} \times \myaxix{\tau(\zeta)}{y}_2,
\{ q \} \times \myaxix{\tau(\zeta)}{y}_3 \big),
\end{align*}
thus $\myaxix{\tau(\xi)}{x} = \myaxix{\tau(\zeta)}{y}$. And in the case $\myaxix{\xi}{x} <_+ \myaxix{\zeta}{y}$, we get again with Corollary \ref{corr_cqxi}
\begin{align*}
( q, \tau(\xi)(x) ) & 
= \myarxix{\Phi^{n(\xi)-1}(\xi)}{x}_1 
\in \myarxix{\Phi^{n(\zeta)-1}(\zeta)}{y}_2
= \{ q \} \times \myaxix{\tau(\zeta)}{y}_2,
\end{align*}
thus $\tau(\xi)(x) \in \myaxix{\tau(\zeta)}{y}_2$. In the same way we see $\tau(\zeta)(y) \in \myaxix{\tau(\xi)}{x}_3$, and $\myaxix{\tau(\xi)}{x} <_+ \myaxix{\tau(\zeta)}{y}$ is proven. $\tau$ is thus an I-scheme, and it is strong according to Theorem \ref{theo_cancel_prod_G}.

\EP

\end{document}